\title{Chiral covers of hypermaps}
\author{Gareth A. Jones\\
School of Mathematics\\
University of Southampton\\
Southampton SO17  1BJ, U.K.\\
{\tt G.A.Jones@maths.soton.ac.uk}
}

\documentclass[10pt]{article}
\usepackage{color}
\usepackage[latin1]{inputenc}
\usepackage{latexsym}
\usepackage{amsmath}
\usepackage{amssymb}
\usepackage{amsfonts}
\newtheorem{thm}{Theorem}[section]

\newtheorem{cor}[thm]{Corollary}

\date{}
\begin{document}

\maketitle

\begin{abstract}
Generalising a conjecture of Singerman, it is shown that there exist orientably regular chiral hypermaps of every non-spherical type. The proof uses the representation theory of automorphism groups acting on homology and on various spaces of differentials. Some examples are given.
\end{abstract}



\section{Introduction}

In 1992, in an unpublished preprint~\cite{Sing}, Singerman conjectured that if $\frac{1}{m}+\frac{1}{n}\le\frac{1}{2}$ there is an orientably regular chiral map of type $\{m,n\}$. Huc\'\i kov\'a, Nedela and \v Sir\'a\v n have announced a proof of this in~\cite{HNS} (see also~\cite[\S3.4]{Sir}). The aim of this note is to prove a similar but more general result for hypermaps of non-spherical type. Whereas most constructions of maps or hypermaps involve group theory or combinatorics, this construction is mainly based on analysis (specifically, spaces of differentials on Riemann surfaces) and representation theory (the action of automorphism groups on such spaces and on associated homology groups). For background on the first topic see~\cite{FK, Jos}, and for the second see~\cite{CR}. It is hoped that such techniques may find further applications in this area.

The author is grateful to David Singerman and to Jozef \v Sir\'a\v n for very helpful comments on earlier drafts of this paper.

\section{Chiral covers}

\begin{thm}\label{mainthm}
If $\mathcal H$ is a finite orientably regular hypermap of non-spherical type, then $\mathcal H$ is covered by infinitely many finite orientably regular chiral hypermaps of that same type.
\end{thm}

\noindent{\sl Proof.} Let $\Delta$ and $\Delta^*$ be the triangle group and extended triangle group of the given type $(l,m,n)$, so $\mathcal H$ corresponds to a torsion-free normal subgroup $N$ of finite index in $\Delta$ with $G:={\rm Aut}^+{\mathcal H}\cong \Delta/N$. This group $N$ is an orientable surface group of genus $g>0$.
Then $\mathcal H$ is reflexible (i.e.~a regular hypermap) if and only if $N$ is normal in $\Delta^*$, in which case $A:={\rm Aut}\,{\mathcal H}\cong\Delta^*/N$. 
\medskip

\noindent{\bf Case 1.} Suppose that $\mathcal H$ is chiral, so $N$ is not normal in $\Delta^*$. Choose any prime $p$ not dividing $|G|$: by a theorem in Euclid's {\sl Elements}, there are infinitely many such primes. Let $N_p=N'N^p$, a torsion-free subgroup which is characteristic in $N$ and hence normal in $\Delta$, and let ${\mathcal H}_p$ be the corresponding finite orientably regular hypermap of type $(l,m,n)$; it covers $\mathcal H$ since $N_p\le N$. Now $|N:N_p|=p^{2g}$, so $N/N_p$ is a normal, and hence characteristic, Sylow $p$-subgroup of $\Delta/N_p$; if $N_p$ were normal in $\Delta^*$ then $N/N_p$ would be normal in $\Delta^*/N_p$, and hence $N$ would be normal in $\Delta^*$, against our assumption. Thus $N_p$ is not normal in $\Delta^*$, so ${\mathcal H}_p$ is chiral. It has genus $p^{2g}(g-1)+1$.

\smallskip

\noindent{\bf Case 2.} Suppose that $\mathcal H$ is reflexible, so that $N$ is normal in $\Delta^*$. Let $e$ be the exponent of $A$. By a theorem of Dirichlet, there are infinitely many primes $p\equiv 1$ mod~$(e)$. Since $p$ is coprime to $|A|$, Maschke's Theorem holds for representations of $A$ over ${\mathbb F}_p$; moreover, since ${\mathbb F}_p$ contains a full set of $e$-th roots of unity,  ${\mathbb F}_p$ is a splitting field for $A$~\cite[Corollary 70.24]{CR}, so the representation theory of $A$ over ${\mathbb F}_p$ is `the same' as that over $\mathbb C$, and this also applies to any subgroup of $A$.

Let $\mathcal S$ be the compact Riemann surface canonically associated with $\mathcal H$, so that the elements of $G$ and of $A\setminus G$ induce conformal and anticonformal automorphisms of $\mathcal S$. The harmonic differentials on $\mathcal S$ form a complex vector space $V=H^1({\mathcal S},{\mathbb C})$ of dimension $2g$ affording a representation $\rho$ of $A$. This space admits a $G$-invariant direct sum decomposition $V^{+}\oplus V^{-}$, where $V^{+}$ and $V^{-}$ are the $g$-dimensional spaces of holomorphic and antiholomorphic differentials on $\mathcal S$, affording complex conjugate representations $\rho^{+}$ and $\rho^{-}$ of $G$. The elements of $A\setminus G$ transpose holomorphic and antiholomorphic differentials, and hence transpose the subspaces $V^{+}$ and $V^{-}$ of $V$.

Integration around closed paths allows one to identify the first homology group  
$H_1({\mathcal S},{\mathbb C})=H_1({\mathcal S},{\mathbb Z})\otimes_{\mathbb Z}{\mathbb C}$
of $\mathcal S$ with the dual space $V^*$ of $V$, so it affords the dual (or contragredient) representation $\rho^*$ of $A$, which is equivalent to the complex conjugate representation $\overline\rho$ since $A$ is finite. Since $\rho\,|_G$ is the sum of two complex conjugate representations, $\overline\rho\,|_G$ is equivalent to $\rho\,|_G$. Thus $H_1({\mathcal S},{\mathbb C})$ also decomposes as a direct sum $H_1^{+}\oplus H_1^{-}$ of two $g$-dimensional $G$-invariant subspaces affording complex conjugate representations $\rho^{+}$ and $\rho^{-}$ of $G$, and these are transposed by elements of $A\setminus G$.

For any prime $p\equiv 1$ mod~$(e)$ let $N_p:=N'N^p$, a characteristic subgroup of $N$ and hence a normal subgroup of $\Delta^*$, and let $M_p:=N/N_p$. Since $N\cong\pi_1{\mathcal S}$ we have $N/N'\cong H_1({\mathcal S},{\mathbb Z})\cong {\mathbb Z}^{2g}$ and hence
\[M_p\cong H_1({\mathcal S},{\mathbb Z})/pH_1({\mathcal S},{\mathbb Z})\cong H_1({\mathcal S},{\mathbb Z})\otimes_{\mathbb Z}{\mathbb F}_p\cong H_1({\mathcal S},{\mathbb F}_p)\cong {\mathbb F}_p^{2g}.\] 
Indeed, these are isomorphisms of ${\mathbb F}_pA$-modules, with the natural action of $A$ on homology corresponding to the induced action of $\Delta^*/N$ by conjugation on $M_p$. By the `isomorphism' between the representation theories of $A$ and its subgroups over $\mathbb C$ and over ${\mathbb F}_p$, the representation $\rho_p$ of $A$ on $M_p$ can be regarded as the reduction mod~$(p)$ of its representation $\overline\rho\sim\rho$ on $V^*=H_1({\mathcal S},{\mathbb C})$, with respect to a suitable basis for $V^*$. In particular, $M_p$ has a $G$-invariant decomposition $M_p^{+}\oplus M_p^{-}$, with $M_p^{+}$ and $M_p^{-}$ affording $g$-dimensional representations of $G$, and the elements of $A\setminus G$ transpose the two direct factors.

The inverse images $N_p^{+}$ and $N_p^{-}$ of $M_p^{+}$ and $M_p^{-}$ in $N$ are torsion-free normal subgroups of finite index $p^g$ in $\Delta$, so let ${\mathcal H}_p^{+}$ and ${\mathcal H}_p^{-}$ be the corresponding finite orientably regular hypermaps of type $(l,m,n)$. These cover $\mathcal H$, and are non-isomorphic as oriented hypermaps since $N_p^{+}\ne N_p^{-}$. Elements of $\Delta^*\setminus G$ transpose $N_p^{+}$ and $N_p^{-}$ by conjugation, so ${\mathcal H}_p^{+}$ and ${\mathcal H}_p^{-}$ form a chiral pair. They have genus $(g-1)p^g+1$. \hfill$\square$



\begin{cor}\label{corollary}
There exist infinitely many orientably regular chiral hypermaps of each non-spherical type.
\end{cor}

\noindent{\sl Proof.} Being residually finite, the triangle group of the given type has a normal subgroup of finite index which contains no non-identity powers of the canonical generators, and is therefore torsion-free. Applying Theorem~\ref{mainthm} to the corresponding orientably regular hypermap gives the required chiral hypermaps. \hfill$\square$

\medskip

\noindent {\bf Comments. 1.} In many cases, the condition $p\equiv 1$ mod~$(e)$ in case~2 is unduly restrictive. It guarantees that every representation of $A$ or $G$ over $\mathbb C$ decomposes in the same way when regarded as a representation over ${\mathbb F}_p$, whereas we are interested in just one representation, namely $\rho$. There are cases, illustrated in the following examples, where $\rho$ decomposes in the required way, yielding a chiral pair, for certain primes $p\not\equiv 1$ mod~$(e)$.

\medskip

\noindent{\bf 2.} The proofs of Theorem~\ref{mainthm} and Corollary~\ref{corollary} make no essential use of the fact that $\Delta$ and $\Delta^*$ are triangle groups, or equivalently that $\mathcal S$ is a Bely\u\i\/ surface, so they yield more general results concerning coverings of compact Riemann surfaces. It is hoped to explore these in a later paper.

\medskip

\noindent{\bf 3.} The method used in case~1 of the proof of Theorem~\ref{mainthm} is sometimes called the `Macbeath trick', since it was used by Macbeath~\cite{Macb} to produce an infinite sequence of Hurwitz groups of the form $\Delta/N'N^m$, with $\Delta$ of type $(2,3,7)$ and $G=\Delta/N\cong L_2(7)$. (See \S4 for background, and for more details concerning this example.) The resulting hypermaps (maps of type $\{3,7\}$) are all regular, whereas here, in cases~1 and 2 of Theorem~\ref{mainthm}, Macbeath's method is modified to produce chiral hypermaps.

\section{An example in genus $2$}

There is a unique orientably regular hypermap $\mathcal H$ of genus $2$ and type $(8,2,8)$. This is a map of type $\{8,8\}$, arising from an epimorphism $\Delta=\Delta(8,2,8)\to G = C_8$ with kernel $N=\Delta'$. By its uniqueness $\mathcal H$ is reflexible, with $A\cong\Delta^*/N\cong D_8$; it appears as R$2.6$ in Conder's list of regular maps~\cite{Con}.

One can construct $\mathcal H$ as a map by identifying opposite sides of an octagon, so that $\mathcal H$ has one vertex, four edges and one face. Going around the boundary of the octagon gives a presentation
\[N=\langle a, b, c, d\mid abcda^{-1}b^{-1}c^{-1}d^{-1}=1\rangle.\]
The images of $a, b, c, d$ in $N/N'$ form a basis for the homology group $V^*=H_1({\mathcal S},{\mathbb C})$ of the underlying surface $\mathcal S$. The automorphism group $A$ of $\mathcal H$  is induced by the isometries of the octagon: thus $G=\langle r\rangle$ and $A=\langle r, s\rangle$ where $r$ is a rotation through $\pi/4$ and $s$ is a reflection; these are represented on homology by the matrices
\[
\left(\,\begin{matrix}&1&&\cr &&1&\cr &&&1\cr -1&&&\cr\end{matrix}\,\right)
\quad{\rm and}\quad
\left(\,\begin{matrix}&&&1\cr &&1&\cr &1&&\cr 1&&&\cr\end{matrix}\,\right).
\]

The Riemann surface $\mathcal S$ underlying $\mathcal H$ is the hyperelliptic curve
\[w^2=z(z^4-1),\]
with Bely\u\i\/ function $\beta:(w,z)\mapsto z^4$. The single vertex is at $(0,0)$, the face-centre is at $(\infty,\infty)$, and the edge-centres are at the four points $(0,z)$ with $z^4=1$. The automorphisms $r$ and $s$ send points $(w,z)\in{\mathcal S}$ to $(\zeta w,iz)$ and $(\overline w,\overline z)$ respectively, where $\zeta=\exp(2\pi i/8)$.

The differentials
\[\omega_1=\frac{dz}{w}\quad{\rm and}\quad\omega_2=\frac{z\,dz}{w}\]
form a basis for $V^{+}$~\cite[\S III.7.5]{FK}, and their conjugates form a basis for $V^{-}$. The rotation $r$ sends $\omega_1$ and $\omega_2$ to $i\omega_1/\zeta=\zeta\omega_1$ and $i^2\omega_2/\zeta=\zeta^3\omega_2$, so these span $1$-dimensional $G$-invariant subspaces $E_{\lambda}$ on which $r$ has eigenvalues $\lambda=\zeta$ and $\zeta^3$; there is a similar decomposition for $V^{-}$, except that here the eigenvalues are $\zeta^{-1}$ and $\zeta^{-3}$. The action of $s$ is to transpose $E_{\zeta}$ with $E_{\zeta^{-1}}$, and $E_{\zeta^3}$ with $E_{\zeta^{-3}}$.

Taking the dual space and then reducing mod~$(p)$, we see that the same applies to the $A$-module $M_p=H_1({\mathcal S},{\mathbb F}_p)$ for any prime $p\equiv 1$ mod~$(8)$, with $\zeta$ now interpreted as a primitive $8$th root of $1$ in the splitting field ${\mathbb F}_p$. This gives a decomposition $M_p=M_p^{+}\oplus M_p^{-}$ with $G$-submodules $M_p^{+}=E_{\zeta}\oplus E_{\zeta^3}$ and $M_p^{-}=E_{\zeta^{-1}}\oplus E_{\zeta^{-3}}$ transposed by $s$. Lifting these back to $N$ gives subgroups $N_p^{+}$ and $N_p^{-}$ corresponding to a chiral pair of maps ${\mathcal H}_p^{+}$ and ${\mathcal H}_p^{-}$ of type $\{8,8\}$ and genus $p^2+1$, as in case~2 of the proof of Theorem~\ref{mainthm}.

\medskip

\noindent{\bf Comments. 1.} In this example, $A$ has exponent $8$, so the smallest prime for which this construction applies is $p=17$, giving a chiral pair of orientably regular hypermaps of type $\{8,8\}$ and genus $290$. They correspond to entry C$290.4$ in Conder's list of chiral maps~\cite{Con}. 

\medskip

\noindent{\bf 2.} This particular example also yields two chiral pairs of $p$-sheeted coverings of $\mathcal H$, corresponding to the four maximal submodules of $M_p$, each omitting one of the four $1$-dimensional eigenspaces $E_{\lambda}$. These are orientably regular maps of type $\{8,8\}$ and genus $p+1$, and each pair consists of the duals of the other; when $p=17$ they correspond to entry C$18.1$ in~\cite{Con}.

\medskip

\noindent{\bf 3.} It is interesting to see what happens if we use primes $p\not\equiv 1$ mod~$(8)$ in this example, so that ${\mathbb F}_p$ is not a splitting field for $A$. If $p\equiv 3$ mod~$(8)$ then $M_p$ is a direct sum of two irreducible $G$-submodules, with $r$ having eigenvalues $\zeta$ and $\zeta^3$ on one, and their inverses on the other, where $\zeta$ is now a primitive $8$th root of $1$ in ${\mathbb F}_{p^2}$; these submodules are transposed by $s$, so we obtain a chiral pair of self-dual maps of type $\{8,8\}$ and genus $p^2+1$; see C$10.3$ in~\cite{Con} for the case $p=3$, and C$122.7$ for $p=11$.

The same applies if $p\equiv -3$ mod~$(8)$, except that $r$ now has eigenvalues $\zeta$ and $\zeta^{-3}$ on one submodule, and their inverses on the other;  again these submodules are transposed by $s$, so we obtain a chiral pair of maps of type $\{8,8\}$ and genus $p^2+1$, now duals of each other; see C$26.1$ and C$170.7$ in~\cite{Con} for $p=5$ and $13$.

If $p\equiv -1$ mod~$(8)$ then $M_p$ is a direct sum of two irreducible $G$-submodules, with $r$ having eigenvalues $\zeta^{\pm 1}$ on one and $\zeta^{\pm 3}$ on the other; these submodules are both invariant under $s$, so we obtain a dual pair of regular maps of type $\{8,8\}$, rather than a chiral pair; see R$50.7$ in~\cite{Con} for the case $p=7$.

Finally, if $p=2$ we again obtain no chiral maps, but there is a unique series
\[M_2>(r-1)M_2>(r-1)^2M_2>(r-1)^3M_2>(r-1)^4M_2=0\]
 of $A$-submodules, giving regular maps of type $\{8,8\}$ and genus $2, 3, 5, 9$ and $17$.

\medskip

\noindent{\bf 4.} Further examples of chiral and regular hypermaps, arising as elementary abelian coverings of regular hypermaps of genus $2$, have been found by Kazaz in~\cite{Kaz}; this example is based on his methods.

\section{Chiral covers of Klein's quartic curve}

Klein's quartic curve $x^3y+y^3z+z^3x=0$ is a compact Riemann surface $\mathcal S$ of genus $3$, which carries a regular map $\mathcal H$ of type $\{3,7\}$ with $G\cong L_2(7)$ and $A\cong PGL_2(7)$; see~\cite{Lev} for a comprehensive study of this curve. The method of proof of case~2 of Theorem~\ref{mainthm}, with $\Delta$ of type $(2,3,7)$, yields chiral pairs of maps of type $\{3,7\}$ and genus $2p^3+1$ for all primes $p\equiv 1$ mod~$(168)$. These are of interest because their automorphism groups, as finite quotients of $\Delta$, are all Hurwitz groups, attaining Hurwitz's upper bound of $84(g-1)$ for the number of automorphisms of a compact Riemann surface of a given genus $g>1$.

In fact, such chiral pairs exist for all primes $p\equiv 1, 2$ or $4$ mod~$(7)$. These can be obtained from a classification by Cohen~\cite{Coh} of those Hurwitz groups which arise as abelian coverings of $G$; his methods of construction are purely algebraic, using $6\times 6$ matrices which can be interpreted as representing generators of $\Delta$ on various homology modules. (These covers were also obtained in an earlier paper of Wohlfahrt~\cite{Woh}, using ideals in the ring of integers of ${\mathbb Q}(\sqrt{-7})$.) The case $p=2$ is due to Sinkov~\cite{Sink}, giving the chiral pair C17.1 of genus $17$ in~\cite{Con}.

In this example, $V^{+}$ and $V^{-}$ are irreducible $G$-submodules of $V$, affording complex conjugate representations with characters taking the values $3, -1, 0, 1$ on elements of orders $1, 2, 3, 4$, and $\alpha=(-1+\sqrt{-7})/2$ and $\overline\alpha$ on the two classes of elements of order $7$. Elkies has studied these representations, and their reduction modulo various primes, in~\cite[\S 1]{Elk}. They give irreducible representations over ${\mathbb F}_p$ for any prime $p$ such that $-7$ is a quadratic residue, that is, for $p\equiv 1, 2$ or $4$ mod~$(7)$; for such primes we have a $G$-module decomposition $M_p=M_p^{+}\oplus M_p^{-}$, giving a chiral pair ${\mathcal H}_p^{+}$ and ${\mathcal H}_p^{-}$. However, no other primes give chiral pairs. The module $M_7$ is indecomposable but reducible, with a submodule of dimension $3$ yielding a regular map of genus $687$, and the zero submodule yielding one of genus $235299$. For primes $p\equiv 3, 5$ or $6$ mod~$(7)$, $M_p$ is irreducible and again no chiral coverings arise. For instance, $M_3$ yields only a regular map of genus $1459$, corresponding to its zero submodule.

\medskip


\end{document}